\documentclass[journal,twoside,web]{ieeecolor}

\usepackage{generic}
\usepackage{cite}
\usepackage{amsmath,amssymb,amsfonts}
\usepackage{algorithmic}
\usepackage{graphicx}
\usepackage{textcomp}
\usepackage{cleveref}

\usepackage{graphicx}
\usepackage{subfigure}
\usepackage{epsfig} % for postscript graphics files
\usepackage{cite}
\usepackage{lcsys}

\newtheorem{theorem}{Theorem}[section]
\newtheorem{lemma}[theorem]{Lemma}

\begin{document}

\def\BibTeX{{\rm B\kern-.05em{\sc i\kern-.025em b}\kern-.08em
    T\kern-.1667em\lower.7ex\hbox{E}\kern-.125emX}}
\markboth{\journalname, VOL. XX, NO. XX, XXXX 2017}
{Author \MakeLowercase{\textit{et al.}}: Preparation of Papers for IEEE Control Systems Letters (August 2022)}

\title{Safe Navigation in the Presence of Range-Limited Pursuers}

\author{
Thomas~Chapman$^{a}$,
Alexander~Von~Moll$^{a}$,
and~Isaac~E.~Weintraub$^{a}$
\thanks{$^{a}$ Air Force Research Laboratory, Aerospace Systems Directorate, Wright-Patterson AFB, OH 45433, \{\texttt{thomas.chapman.6, alexander.von\_moll, isaac.weintraub.1}\} \texttt{@us.af.mil}}
\thanks{This paper is funded in part by AFOSR LRIR 24RQCOR002. DISTRIBUTION STATEMENT A. Approved for public release. Distribution is unlimited. AFRL-2025-4392; Cleared 03-Sep-2025.}
}

\maketitle
\thispagestyle{empty}

%%%%%%%%%%%%%%%%%%%%%%%%%%%%%%%%%%%%%%%%%%%%%%%%%%%%%%%%%%%%%%%%%%%%%%%%%%%%%%%%
\begin{abstract}
    This paper examines the degree to which an evader seeking a safe and efficient path to a target location can benefit from increasing levels of knowledge regarding one or more range-limited pursuers seeking to intercept it. Unlike previous work, this research considers the time of flight of the pursuers actively attempting interception. It is shown that additional knowledge allows the evader to safely steer closer to the threats, shortening paths without accepting additional risk of capture. A control heuristic is presented, suitable for real-time implementation, which capitalizes on all knowledge available to the evader.
\end{abstract}


\section{Introduction}
\IEEEPARstart{N}{avigation} in an obstacle-laden environment represents one of the most consequential problems for autonomous systems, particularly when those obstacles are nonstationary or nondeterministic. This fundamental challenge has garnered significant attention across multiple domains, from military applications involving unmanned aerial vehicles to civilian scenarios such as collision avoidance in congested airspace \cite{faa2023uam}.

\subsection{Prior Navigation Work}
Significant effort has been devoted to developing robust methods of navigation that intelligently avoid obstacles. These obstacles may be static \cite{chu2012local}, dynamic \cite{saunders2005static}, or uncertain \cite{stagg2025probabilistic, tordesillas2019faster, costa2019survey}.

However, most of these techniques assume that the dynamical obstacles move without regard to the navigating agent. One conservative way to address this is to model obstacles as mobile adversaries engaged in pursuit of the navigating agent. Safe navigation in such scenarios requires understanding not just where threats might move, but how they will deliberately maneuver to achieve interception. 

\subsection{Pursuit-Evasion Differential Games}
Differential game theory provides a mathematical framework for analyzing strategies between competing agents \cite{isaacs1965differential}.
Pursuit-evasion dynamics have been studied extensively in many contexts including those with imperfect state information \cite{rhodes1969differential}. 

A key insight that has emerged is the importance of information to an agent's strategy. Pursuit-evasion games with perfect information \cite{chaudhari2023time, morgan2016interception} operate very differently from those with imperfect or asymmetric information \cite{bopardikar2008discrete}. By exploiting information limitations, agents are able to induce suboptimal behavior from their adversaries \cite{shishika2024deception}, which can significantly alter their optimal strategy.  

\subsection{Reachability Analysis}
The domain that unites navigation and pursuit-evasion problems is reachability analysis, which provides a spatial framework for understanding both passive and adversarial threat avoidance.
In navigation contexts, it has been used to extend the concept of velocity obstacles \cite{fiorini1998motion} to handle unpredictable, dynamically constrained obstacles \cite{wu2012guaranteed}. 

When applied to pursuit-evasion scenarios, reachability analysis describes not just where an obstacle might passively drift, but where an active pursuer can deliberately position itself to intercept. This has been examined for the case of range-limited pursuers with and without a capture radius \cite{weintraub2023range}.
A review of reach-avoid differential games with simple motion is provided in \cite{yan2023multiplayer}.

\subsection{The Engagement Zone}
The concept of an engagement zone (EZ) provides a systematic approach for determining when and where interception is possible between range-limited pursuing threats and navigating agents \cite{vonmoll2024basic}.
This has been leveraged in practical navigation scenarios involving adversarial pursuers \cite{weintraub2022optimal, dillon2023optimal}.

Importantly, these works have focused only on cases where the pursuers represent potential, rather than active, threats. That is, the evader navigates based on the theoretical danger a pursuer would pose if it chose to engage, rather than responding to pursuers that are actively attempting to intercept. 

The present work extends the framework of EZ avoidance to scenarios where pursuers actively engage, and describes how an evader can exploit knowledge of its adversaries to achieve more efficient navigation while maintaining safety guarantees. 

\subsection{Organization}
This paper is organized as follows: In \Cref{sec-concepts}, the knowledge levels, agent models, and maximum cutting angle for safe vehicle planning are defined and described. The safe headings for each knowledge level are explained in \Cref{sec-kl1}, \Cref{sec-kl2}, and \Cref{sec-kn3}. In \Cref{sec-path-planning}, a heuristic for evading capture is formulated and its value is explained. A simulation is conducted in \Cref{sec-sim} with results showing the improvement made to the path plan when more knowledge is available. Lastly, final conclusions and prospects for future work are provided in \Cref{sec-conclusions}.

\subsection{Contributions}
This work examines navigation of an evader in the presence of one or more adversarial pursuers. It analyzes three levels of knowledge an evader might reasonably have about its pursuer(s), and discusses the means by which the evader can take advantage of this knowledge to reduce its path length while remaining safe from interception. A control heuristic is developed which is suitable for real time implementation. An example navigation scenario is included and discussed. This paper builds upon the notion of the engagement zone \cite{vonmoll2024basic} to encompass actively maneuvering pursuers; also, it details how safe trajectories change as the levels of information about pursuers change.

\section{Concepts}
\label{sec-concepts}

\subsection{Knowledge in Path Planning}
The ability of an evader to safely and efficiently navigate to a target location in the presence of pursuers is intimately related to the knowledge it possesses about said pursuers. For instance, suppose the evader only knows the range and initial location of the pursuers. In this case, the only way to \textit{ensure} the safety of its route is to never enter any pursuer's reachability region (RR) . The RR is defined as the area of space where a pursuer can operate or be effective. [Safety is defined as evading capture from incoming pursuers. Point capture is considered and there is no modeled probability of capture.

A pursuer is said to have `launched' when it has begun moving from its initial location in an attempt to intercept the evader. If the evader is aware that a pursuer has not yet launched,b it is capable of adopting trajectories that travel through that pursuer's RR which cannot be intercepted by any subsequent launch (the pursuer will be `outrun'). 
These trajectories can be determined based on those headings that remain outside the engagement zone of the pursuer. 

This work examines three possible levels of knowledge the evader could possess:

\textbf{Level 1}: The evader is only aware of its pursuers' initial locations and motion parameters.

\textbf{Level 2}: The evader is additionally aware of the instant(s) a pursuer chooses to launch.

\textbf{Level 3}: The evader is additionally aware of the positions of all launched pursuers. 

These levels are cumulative, with each subsequent level containing all information in the previous.

\subsection{Pursuer and Evader Models}

The evader and pursuer(s) each move in the 2D-plane at constant speeds $v_E$ and $v_P$ respectively and can change their headings instantaneously. Let $E = [x_E,y_E] \in \mathbb{R}^2$ and $P = [x_P,y_P] \in \mathbb{R}^2$ be the states of the evader and pursuer. The dynamic equations of motion for the evader and pursuer are respectively $[\dot{x}_E,\dot{y}_E] = v_E [\cos \psi, \sin \psi]$, and $[\dot{x}_P,\dot{y}_P] = \vec{v}_P$ where $||\vec{v}_P|| = v_P$. Scenarios are characterized by the speed ratio $\mu = v_E/v_P$. For the entirety of this work, it is assumed that $\mu \leq 1$. The pursuer(s) are time-limited (akin to a fuel limitation) to a maximum $t_\mathrm{max}$ from their moment of launch. This amounts to a maximum range of $R = v_P t_{\mathrm{max}}$. Thus, a pursuer's reachability region is a circle of radius $R$ centered at its initial location. For each scenario examined in this work, all pursuers are initially located at the same location.

\subsection{The Engagement Zone and Maximum Cutting Angle}

\begin{figure}[htbp]
    \centering
    \includegraphics[width=1.0\linewidth]{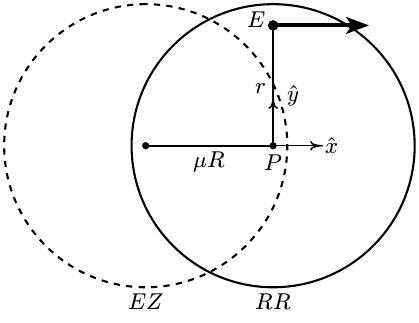}
    \caption{An engagement zone for a simple motion pursuer where $\mu < 1$ is shown. The pursuer and evader are located at points $P$ and $E$ respectively. The evader's velocity in the $+\hat{x}$ direction is shown, and the EZ is consequently shifted in the $-\hat{x}$ direction. Because the evader is outside the EZ, it will be safe if it maintains its current heading, even if the pursuer launches. The evader will outrun the pursuer in this case.}
    \label{fig-ez}
\end{figure}

The engagement zone (EZ) of a pursuer determines which headings can be safely adopted while within its RR. When a pursuer is faster or equal in speed to an evader (i.e. $\mu \leq 1$), the EZ is simply the circle of the pursuer's RR, shifted a distance $\mu R$ in the direction opposite the evader's heading \cite{vonmoll2024basic}. By definition, if the evader is outside a pursuer's EZ, it can safely continue along its current heading indefinitely without being intercepted (provided it remains outside the EZ along this trajectory). This can be seen in \Cref{fig-ez}.

\begin{figure}[htbp]
    \centering
    \includegraphics[width=1.0\linewidth]{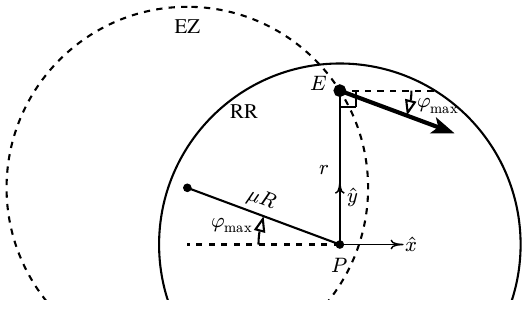}
    \caption{An evader taking the maximum cutting angle $\varphi_\mathrm{max}$ for a pursuer where $\mu < 1$ is shown. By construction, this heading results in the evader lying on the  boundary of the engagement zone.}
    \label{fig-max-cutting-angle}
\end{figure}

It follows naturally to ask what headings an evader at a given location can take such that it will not enter a pursuer's EZ. Consider \Cref{fig-max-cutting-angle}.

\begin{lemma}
    The maximum heading change $\varphi_\mathrm{max}$ (w.r.t a heading that is normal to the line between pursuer and evader) that the evader can take before ending inside said pursuer's EZ is
    \begin{equation}
    \label{eq-phi_max-general}
        \varphi_\mathrm{max} = \mathrm{arcsin} \left(\frac{(\mu^2 - 1)R^2 + r^2}{2 \mu r R}\right).
    \end{equation}
\end{lemma}

\begin{proof}
    Due to the rotational symmetry present, the evader's location can be taken as $(0,r)$ without loss of generality (see \Cref{fig-max-cutting-angle}). By definition, a heading change of $\varphi_\mathrm{max}$ will correspond to the evader being on the boundary of the EZ. The equation of the EZ is given by a circle shifted a distance $\mu R$ from the RR center in the direction opposite the evader's heading \cite{vonmoll2024basic}:  

\begin{equation}
\label{eq-EZ-circle-equation}
    (x + \mu R \cos \varphi_\mathrm{max})^2 + ( y - \mu R \sin \varphi_\mathrm{max})^2 = R^2.
\end{equation}

Thus, \Cref{eq-phi_max-general} follows by substituting the evader's location into \Cref{eq-EZ-circle-equation} and solving for $\varphi_\mathrm{max}$.

\end{proof}

In the case where the evader lies on the boundary of the RR (i.e. $r = R$), \Cref{eq-phi_max-general} simplifies to: 
\begin{equation}
\label{eq-phi_max-RR-boundary}
    \varphi_\mathrm{max} = \mathrm{arcsin}\left( \frac{\mu}{2}\right).
\end{equation}

It is important to state that $\varphi_\mathrm{max}$ is measured relative to the local perpendicular of the vector from the pursuer's location to the evader's location. This must be accounted for to accurately track the evader's safe headings relative to the fixed reference of the $\hat{x}$ and $\hat{y}$ axes. 

To distinguish it from the absolute heading $\psi$, the variable $\varphi$ is defined as the \emph{cutting angle} because it describes how much the evader can safely `cut' into a pursuer's RR. A positive value indicates than the evader can safely approach a pursuer, and a negative value implies that the evader must turn away to remain safe. If the argument of $\mathrm{arcsin}(\cdot)$ in \Cref{eq-phi_max-general} is outside the domain $[-1,1]$, no heading angle exists that will result in the evader being on the boundary of the EZ. This is possible when $r$ is very large (the evader is outside the EZ for all headings) or very small (the evader is inside the EZ for all headings).

Note that, if the pursuer launches on collision course at the same instant the evader in \Cref{fig-max-cutting-angle} begins moving in the direction indicated by $\varphi_\mathrm{max}$, it would intercept the evader on the boundary of the RR (i.e. a max-range interception) because the evader begins on the boundary of the EZ. For the purposes of this work, such interceptions are neglected as they can be theoretically prevented by infinitesimally overestimating the pursuer's time limit. 

\section{Knowledge Level 1}
\label{sec-kl1}
The lowest level of knowledge considered is when the evader is only aware of the initial locations and motion parameters of its pursuers. This is the absolute minimum amount of information necessary to develop a safe path plan.

In this case, the evader must never enter the RRs of any of its pursuers. This is because it cannot guarantee that it will avoid interception, since pursuers may have launched at some earlier moment and be anywhere within their RRs at any given time. 
The minimum-time paths, subject to this constraint, are the longest compared to subsequent knowledge levels.

\section{Knowledge Level 2}
\label{sec-kl2}
Suppose, in addition to knowledge level 1, the evader also has an indication of launch for each pursuer. That is, the evader knows, for each pursuer, the time $\Delta t$ in the past when they launched. For pursuers who have not yet launched, $\Delta t = 0$.

It is now possible to safely cut into the RR of a pursuer because the evader is capable of bounding the pursuer's possible locations within its RR. The question is how sharp a cut can be made. There are three possibilities to consider:

\subsection{$\Delta t = 0$}

If a pursuer has not yet launched, its location is known. Thus, the maximum cutting angle is given by \Cref{eq-phi_max-general}.

\subsection{$0 < \Delta t < t_\mathrm{max}$}

If a pursuer has launched, the question is more complex. Since, at knowledge level 2, the evader is not aware of the pursuers' in-motion positions, it can only describe where a launched pursuer \emph{could} be: anywhere within a circle of radius $ s = v_P \Delta t$ around its initial location. The evader would like to know the largest cutting angle that is safe from a pursuer who is \emph{anywhere} in that region. 

\begin{lemma}
\label{lem-2}
      The maximum cutting angle (measured w.r.t. the normal of the line from the RR center to the evader) that is safe from a pursuer who has launched $\Delta t$ in the past under knowledge level 2 (i.e., no information on the pursuer's current position) is

\begin{equation}
\label{eq-phi_max-nonzero-launch-time}
    \varphi_\mathrm{max}^\mathrm{L2} = \mathrm{arcsin} \left(\frac{\mu^2(R-s)^2 +r^2 - R^2}{2 \mu r (R - s)}\right).
\end{equation}

The superscript L2 in \Cref{eq-phi_max-nonzero-launch-time} corresponds to knowledge level 2.

\end{lemma}

\begin{proof}
    Intuitively, the most challenging pursuer position must be on the boundary of the disk of radius $s$ centered on the pursuer's initial position. Consequently, the pursuer's possible positions can be parameterized with an angle $\gamma$ as shown in \Cref{fig-cutting-angle-kl2}.

\begin{figure}[htbp]
    \centering
    \includegraphics[width=0.9\linewidth]{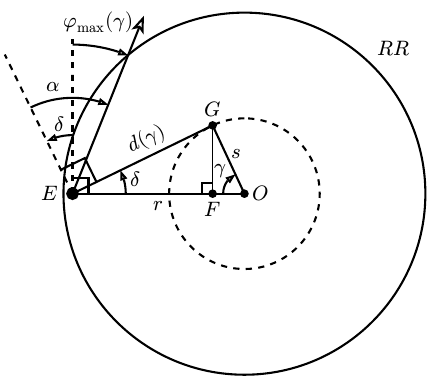}
    \caption{The maximum safe cutting angle for a knowledge level 2 evader given that a pursuer launched $\Delta t$ before the evader reached its current location $E$ is shown. The pursuer could be anywhere in a circle of radius $s = v_P \Delta t$ around its initial location. Its possible positions are parameterized with the angle $\gamma$. Each possible pursuer position has an associated maximum cutting angle $\varphi_\mathrm{max} (\gamma)$.}
    \label{fig-cutting-angle-kl2}
\end{figure}

From there, an equation can be derived for the maximum cutting angle w.r.t. a pursuer at angle $\gamma$:

\begin{equation}
\label{eq-kl2-phi-of-gamma}
    \varphi_\mathrm{max} (\gamma) =  \alpha  - \delta,
\end{equation}
with $\alpha$ and $\delta$ given by

\begin{equation}
\label{eq-kl2-alpha-delta}
    \begin{aligned}
    \alpha &= \mathrm{arcsin}\left(\frac{(\mu^2-1) (R-s)^2 + d(\gamma)^2)}{2 \mu d(\gamma) (R-s)}\right),\\
    \delta &= \mathrm{arctan2}(s \sin \gamma, r - s \cos \gamma),
    \end{aligned}
\end{equation}

where $r$ is the radial distance of the evader from the RR center and $d(\gamma) = \sqrt{(r^2 - 2 r s \cos \gamma + s^2)}$ is the distance between a pursuer at angle $\gamma$ and the evader.
The quantity $\alpha$ is simply the application of \Cref{eq-phi_max-general} to a pursuer at angle $\gamma$ with range $R-s$, and $\delta$ corrects this cutting angle to be measured relative to a common perpendicular (that of the line between RR center and evader position). It is derived from an arc-tangent on the $\triangle
 E F G$ (shown in \Cref{fig-cutting-angle-kl2}). Without the second term, the cutting angles for each $\gamma$ would be measured with respect to a different local normal, and thus would not be comparable.

The minimum of this function can be obtained by differentiating \Cref{eq-kl2-phi-of-gamma} with respect to $\gamma$ and solving for the $\gamma^*$ that satisfies $\varphi'_\mathrm{max} (\gamma^*) = 0$. It is given by:

\begin{equation}
\label{eq-kl2-gamma-critical}
    \gamma^*  = \mathrm{arccos}\left( \frac{(R^2 + r^2 - \mu^2(R-s)^2)}{2 r R}\right).
\end{equation}

Finally, the result follows from substituting $\gamma^*$ into \Cref{eq-kl2-phi-of-gamma}.
\end{proof}

\subsection{$\Delta t \geq t_\mathrm{max}$}
Since the pursuers are time-limited, any that launched more than $t_\mathrm{max}$ in the past can be ignored.

\section{Knowledge Level 3}
\label{sec-kn3}
Suppose the evader additionally has full knowledge of each pursuer's time-varying position after launch. This is the highest level of available information considered. 

\begin{lemma}
\label{lem-1}
    The maximum cutting angle (measured w.r.t the normal of the line from pursuer to evader) that is safe from a pursuer under knowledge level 3 (i.e. the pursuer's position is known) is
\begin{equation}
\label{eq-kl3-phi-max}
    \varphi_\mathrm{max}^\mathrm{L3} = \mathrm{arcsin}\left(\frac{(\mu^2-1)(R-s)^2 + d^2)}{2 \mu d (R-s)}\right).
\end{equation}
As before, the superscript indicates the knowledge level.
\end{lemma}

\begin{proof}
  Knowing the location of the pursuer makes this geometry equivalent to that considered for \Cref{lem-1} (where the pursuer is known to be at its initial location), but with a shifted launch location and reduced pursuer range.
  Thus, $\varphi_\mathrm{max}^\mathrm{L3}$ is given by a modification of \Cref{eq-phi_max-general} where the pursuer range, $R$, has been replaced with $R-s$ and distance from RR center to evader, $r$, is replaced with the distance from pursuer to evader, $d$. Note that $\varphi_\mathrm{max}^\mathrm{L3}$ is measured with respect to the normal of the line between pursuer and evader.
\end{proof}

It is important to recognize that the benefit of this additional knowledge is dependent on pursuers taking suboptimal trajectories. To ensure safety, an evader with knowledge level 2 must assume its pursuers adopt the most challenging paths (in the sense that they minimize the maximum safe cutting angle) regardless of their actual behavior (because this is not known). If the pursuers indeed \emph{do} take these paths, a knowledge level 3 evader will follow exactly the same trajectory as a knowledge level 2 evader, and thus see no additional benefit.

\section{Path Planning Heuristic} 
\label{sec-path-planning}

A heuristic for path planning can now be formulated: to navigate through an RR, move at the maximum cutting angle $\varphi_\mathrm{max}$, unless it is greater than or equal to the cutting angle needed to go straight to the target location, then simply do so.
Note that this heuristic can be followed through either clockwise or counter-clockwise motion about the RR center. A rule to select the turn direction may be defined, but the direction is selected to be clockwise \emph{a priori} for the purposes of this work.

\begin{figure}[htbp]
    \centering
    \includegraphics[width=0.9\linewidth]{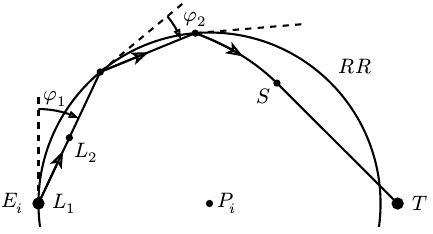}
    \caption{An example path following the presented heuristic is shown. A knowledge level 2 evader starts at $E_i$ and wants to reach the target location $T$. Immediately after starting, a pursuer launches from $P_i$ when the evader is at $L_1 = E_i$. The evader travels straight to the RR at cutting angle $\varphi_1$ given by \Cref{eq-phi_max-RR-boundary}, and a second pursuer launches when it is at $L_2$. Upon reaching the RR, it selects a new cutting angle $\varphi_2$ given by \Cref{eq-phi_max-nonzero-launch-time}, and again proceeds until reaching the boundary. No more pursuers launch, so taking the maximum cutting angle, given by \Cref{eq-phi_max-general}, results in the evader moving with non-constant heading until it reaches $S$ and can proceed straight to the target.}
    \label{fig-heuristic-guide}
\end{figure}

\begin{figure}[htbp]
    \centering
    \includegraphics[width=0.55\linewidth]{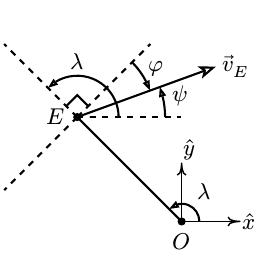}
    \caption{A depiction of the relationship between heading $\psi$, cutting angle $\varphi$, and line-of-sight angle $\lambda$ is shown. In this work, the origin $O$ is the center of the pursuer's RR.}
    \label{fig-heading-and-cutting-angle}
\end{figure}

To understand the implications of this heuristic, consider \Cref{fig-heuristic-guide}. A knowledge level 2 evader starts on the boundary of an RR, and must reach a target location $T$ on the other side. Initially, no pursuer has yet launched. The evader begins by moving at the maximum cutting angle $\varphi_1$, given by \Cref{eq-phi_max-RR-boundary}. A pursuer instantly proceeds to launch, so the evader continues straight until reaching the RR's boundary, traveling along a chord of the circle. This is because, at knowledge level 2, always moving at the maximum cutting angle results in a straight line path when a pursuer launches. To see this, consider \Cref{fig-heading-and-cutting-angle}. 

The absolute heading of an evader can be written as:

\begin{equation}
\label{eq-absolute-heading}
    \psi = \lambda - \varphi - \frac{\pi}{2}. 
\end{equation}

Taking the time derivative of both sides yields:

\begin{equation}
\label{eq-absolute-heading-time-derivative}
    \dot{\psi} = \dot{\lambda} - \dot{\varphi}_\mathrm{max} = -v_E \frac{\cos \varphi_\mathrm{max} }{r} - \dot{\varphi}_\mathrm{max},
\end{equation}

where the chosen cutting angle $\varphi$ is assumed to be the maximum. This can be shown to be zero algebraically when using  $\varphi_\mathrm{max}^\mathrm{L2}$ from \Cref{eq-phi_max-nonzero-launch-time}. A constant heading results in a straight line path.

Moreover, this straight line path will reach the boundary of the RR at the same moment the pursuer (whose true position is unknown at knowledge level 2) \emph{could} have reached the boundary (i.e. when $s = R \Leftrightarrow \Delta t = t_\mathrm{max}$). 
This is implicit in the formulation of the maximum cutting angle. 
By construction, moving at this angle keeps the evader on the EZ boundary of its pursuer. This prevents all but max-range interceptions (which are neglected for reasons discussed previously).

Thus, if there was initially only one pursuer at the RR's center, the evader could safely proceed directly to the target location, as said pursuer would have just exhausted its time limit. If there were multiple pursuers, the evader will simply continue to move at the maximum cutting angle.

If, as in \Cref{fig-heuristic-guide}, a pursuer is still engaged in pursuit when the evader reaches the boundary, the evader will reenter the RR at a new maximum cutting angle $\varphi_2$ given by \Cref{eq-phi_max-nonzero-launch-time}. The evader travels straight until reaching the boundary again, at which point said pursuer has exhausted its time limit. If multiple pursuers are simultaneously engaged in pursuit, each will have an associated maximum cutting angle (along with the theoretical unlaunched pursuer at the origin), and the evader should choose the most limiting of them so as to evade capture by \emph{any} of its pursuers.

If no pursuers launch, as in the latter half of the trajectory depicted in \Cref{fig-heuristic-guide}, by always traveling at the maximum cutting angle the evader will follow a smooth curve through the RR until the maximum cutting angle corresponds to the straight line path to the target (point $S$). From there, the evader does not have to consider any subsequent launches as it will simply outrun them. The final appearance of a trajectory will depend on the behavior of the pursuers. If many launch, it will be a series of chords. If none launch, it will be a smooth curve through the RR.

Without specific rules of pursuer launch and guidance behavior, it is not possible to claim any path as 'optimal.' Any given path plan, with the exception of avoiding the RR entirely, is inherently provisional, as it remains viable only until the pursuer(s) act in a manner that necessitates a deviation. 

The value of the presented heuristic is twofold. First, it results in safe control through a pursuer's RR that can be computed and updated in real time. Path planning methods typically calculate the entire path prior to an agent actually beginning it, and thus are ill-suited to handle deviations that become necessary during execution. Second, it does not place excessive requirements on the evader to actively \emph{evade} incoming threats of capture. The evader's response to any pursuer launch is to either ignore it or to stop turning and continue straight. 

It is also possible to define a more aggressive heuristic which attempts to take cutting angles of greater than $\varphi_\mathrm{max}$, but these by necessity require evasive maneuvering by the evader (changing heading) upon pursuer launch. It is an open question whether such a behavior can reduce path lengths in more realistic engagements.

\section{Simulation}
\label{sec-sim}

\begin{figure}[htbp]
    \centering
    \includegraphics[width=0.9\linewidth]{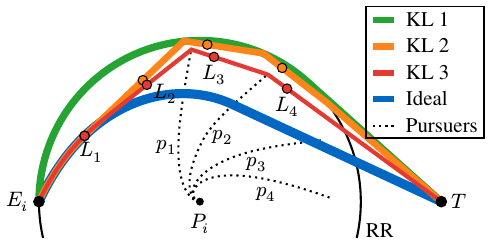}
    \caption{Multiple evader trajectories through a RR following the cutting angle heuristic at different knowledge levels (KL) are shown. The dots on the KL2 and KL3 curves represent the locations of those evaders at the moments a pursuer launches. Those pursuers that seek the KL3 evader have their trajectories labeled $p_1,p_2,p_3,p_4$. These trajectories begin when the KL3 evader is at $L_1,L_2,L_3,L_4$ respectively. Trajectories for the pursuers of the KL2 evader are not shown because this evader is not aware of the locations of its pursuers. The launch indications are omitted for KL1 because this evader doesn't need to react to launches. The KL1 trajectory (green) fully circumnavigates the RR until exiting tangentially. The KL2 trajectory (orange) follows a series of chords after the initial launch. The KL3 trajectory (red) takes advantage of the guidance limitations of the pursuers and can take a tighter path within the RR. The ideal trajectory (blue) assumes no launches occur. Knowledge levels 2 and 3 both follow the ideal path until the first launch forces them to deviate. }
    \label{fig-simulation}
\end{figure}

This section details an identical engagement treated at all three knowledge levels (see \Cref{fig-simulation}). Simulations were performed on an AMD Ryzen 5 Pro 7530U CPU with 64GB RAM.

The pursuers' initial locations are treated as the origin $(0,0)$, with the evader beginning at the point $(-R,0)$ on the boundary of the pursuers' RR. The evader seeks to reach target location $(q,0)$ with $q \geq R$ in minimum time without being intercepted. A new pursuer launches every $t_L$ after the evader enters the RR while it remains inside the RR. Each pursuer is equipped with standard proportional-navigation guidance \cite{shneydor1998missile} with gain $N$, and is initially oriented in the $-\hat{x}$ direction.

These limitations on the pursuer were chosen deliberately to illustrate when knowledge level 3 can achieve significant reductions in path length over knowledge level 2. If the pursuer takes the worst case trajectory (from the evader's perspective), the evader strategies for knowledge levels 2 and 3 are identical. Any variation from the pursuer's worst case trajectory improves the relative performance of an evader with knowledge level 3 over one with knowledge level 2.

For these results, the parameters selected were $v_E=1$, $\mu = 0.85$, $R=1.7$, $q = 2.55$, $t_L = 0.85$, and $N=6$. The path lengths are summarized in \Cref{sim-results}, and the paths are shown in \Cref{fig-simulation}.
The minimum computation time for computing the evader's safe heading was 0.2 $\mu$s.

Knowledge level 2 results in a 3.78\% reduction in path length over knowledge level 1. Knowledge level 3 shows a greater 7.09\% reduction, and the ideal path, resulting from a scenario where no pursuer launches, is 14.36\% shorter. These results demonstrate that additional knowledge of a pursuer's capabilities allows an evader to more fully exploit the pursuer's limitations to plan shorter paths without additional risk.

A larger speed ratio allows evaders to cut more aggressively through the RR, improving knowledge levels 2, 3, and ideal paths.
An increase in $R$ without accompanying increase in $v_E$ results in more launched pursuers, causing knowledge level 2 and 3 paths to more closely approximate the knowledge level 1 path. Fewer launched pursuers cause the knowledge level 2 and 3 paths to more closely approximate the ideal path.

% Speaking generally, a larger speed ratio allows the evader to cut more aggressively into the RR because it decreases the speed advantage of the pursuer. More aggressive pursuer launching (smaller $t_L$) 

\begin{table}[]
    \centering
    \caption{Path lengths for each trajectory}
    \begin{tabular}{l c c}
        \hline
        Trajectory  & Length & \% Improvement \\ \hline
        Knowledge level 1 & $5.812$ & - \\
        Knowledge level 1 & $5.592$ & $3.78$\% \\
        Knowledge level 1 & $5.399$ & $7.09$\% \\
        Ideal (no launch) & $4.977$ & $14.36$ \% \\ \hline
    \end{tabular}
    \label{sim-results}
\end{table}

\section{Conclusions}
\label{sec-conclusions}

This work formulates the concept of a cutting angle that keeps an evader safe from capture while within a pursuer's reachability region. The cutting angle is used to inform a path planning heuristic which is implemented in simulation.
This heuristic allows for shorter paths without increased risk, and becomes increasingly effective when the evader possesses additional knowledge about its threats.

Future work will expand these findings to encompass scenarios where the evader is faster than its pursuer(s) or is turn constrained, and where the pursuer(s) have variable speed. Additionally, the applicability of the heuristic should be assessed for engagements where pursuers have multiple initial locations such that their reachability regions overlap. Formal robustness analysis and expanded simulation studies are identified as future work.

\bibliographystyle{IEEEtran}
\bibliography{Ref}

\end{document}